\documentclass[12pt,reqno,a4wide]{amsart}

\usepackage{tikz} 
\usepackage{calrsfs}
\usepackage{mathrsfs}

\usepackage{array}
\usepackage{hyperref}

\usetikzlibrary{decorations.pathreplacing}
\usetikzlibrary{fit}

\usepackage{pgfplots}

\allowdisplaybreaks

\catcode`,\active

\catcode`\,12

\begin{document}
	\bibliographystyle{plain}
	
	%
	%
	
	\title
	{Skew 2-Dyck paths via the kernel method}

	\author[H. Prodinger ]{Helmut Prodinger }
	\address{Department of Mathematics, University of Stellenbosch 7602, Stellenbosch, South Africa
		and
		NITheCS (National Institute for
		Theoretical and Computational Sciences), South Africa.}
	\email{warrenham33@gmail.com}


	\begin{abstract}
		We continue on a recent concept introduced by Kariuki and Okoth, about skew 2-Dyck paths, introducing an additional down-step $L$, together with the usual steps
		$U$ (up) and $D$ (down). There is the syntactical condition that $UL$ and $LU$ can never occur. An automaton that checks these conditions is introduced, and the
		relevant generating functions are obtained by applying the kernel method to three functional equations. It is briefly discussed how the setting can be extended to $t$-Dyck paths. As a benefit, prefixes of skew $t$-Dyck paths are also enumerated.
		An approach that scans  2-Dyck paths from right to left is also discussed.
	\end{abstract}

	\maketitle

\section{2-Dyck paths}

Dyck paths consist of up-steps $U=(1,1)$ and down-steps $D=(1,-1)$, start at the origin, never go below the $x$-axis and end at the $x$-axis. If a Dyck path does not necessarily end at the $x$-axis we call it a \emph{prefix} of a Dyck path.

A standard reference on the enumeration of lattice paths is \cite{BF}.

The 2-Dyck paths in the section title consist of up-steps $U=(1,1)$ and down-steps $D=(1,-2)$, never go below the $x$-axis and return to the $x$-axis. They will be the main subject of this paper, together with the epitheton \emph{skew}. To clarify, we present a list of all 12 2-Dyck paths of 9 steps (length 9).

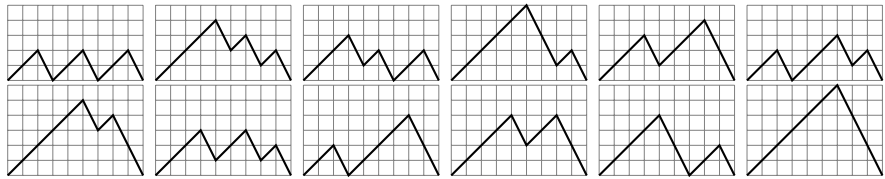
\begin{figure}[h]
	\begin{tikzpicture}[scale=0.2]
		\draw[help lines] (0,0) grid (9,5);
		
		\draw [thick] (0,0) -- (1,1)--(2,2)--(3,0) --(4,1)--(5,2)--(6,0)--(7,1)--(8,2)--(9,0);
		
	\end{tikzpicture}
\begin{tikzpicture}[scale=0.2]
	\draw[help lines] (0,0) grid (9,5);
	
	\draw [thick] (0,0) -- (1,1)--(2,2)--(3,3) --(4,4)--(5,2)--(6,3)--(7,1)--(8,2)--(9,0);

\end{tikzpicture}
\begin{tikzpicture}[scale=0.2]
	\draw[help lines] (0,0) grid (9,5);
	
	\draw [thick] (0,0) -- (1,1)--(2,2)--(3,3) --(4,1)--(5,2)--(6,0)--(7,1)--(8,2)--(9,0);

\end{tikzpicture}
\begin{tikzpicture}[scale=0.2]
	\draw[help lines] (0,0) grid (9,5);
	
	\draw [thick] (0,0) -- (1,1)--(2,2)--(3,3) --(4,4)--(5,5)--(6,3)--(7,1)--(8,2)--(9,0);
	\end{tikzpicture}
\begin{tikzpicture}[scale=0.2]
	\draw[help lines] (0,0) grid (9,5);
	
	\draw [thick] (0,0) -- (1,1)--(2,2)--(3,3) --(4,1)--(5,2)--(6,3)--(7,4)--(8,2)--(9,0);
\end{tikzpicture}
\begin{tikzpicture}[scale=0.2]
	\draw[help lines] (0,0) grid (9,5);
	
	\draw [thick] (0,0) -- (1,1)--(2,2)--(3,0) --(4,1)--(5,2)--(6,3)--(7,1)--(8,2)--(9,0);
\end{tikzpicture}

\begin{tikzpicture}[scale=0.2]
	\draw[help lines] (0,0) grid (9,6);
	
	\draw [thick] (0,0) -- (1,1)--(2,2)--(3,3) --(4,4)--(5,5)--(6,3)--(7,4)--(8,2)--(9,0);
\end{tikzpicture}
\begin{tikzpicture}[scale=0.2]
	\draw[help lines] (0,0) grid (9,6);
	
	\draw [thick] (0,0) -- (1,1)--(2,2)--(3,3) --(4,1)--(5,2)--(6,3)--(7,1)--(8,2)--(9,0);
\end{tikzpicture}
\begin{tikzpicture}[scale=0.2]
	\draw[help lines] (0,0) grid (9,6);
	
	\draw [thick] (0,0) -- (1,1)--(2,2)--(3,0) --(4,1)--(5,2)--(6,3)--(7,4)--(8,2)--(9,0);
\end{tikzpicture}
\begin{tikzpicture}[scale=0.2]
	\draw[help lines] (0,0) grid (9,6);
	
	\draw [thick] (0,0) -- (1,1)--(2,2)--(3,3) --(4,4)--(5,2)--(6,3)--(7,4)--(8,2)--(9,0);
\end{tikzpicture}
\begin{tikzpicture}[scale=0.2]
	\draw[help lines] (0,0) grid (9,6);
	
	\draw [thick] (0,0) -- (1,1)--(2,2)--(3,3) --(4,4)--(5,2)--(6,0)--(7,1)--(8,2)--(9,0);
\end{tikzpicture}
\begin{tikzpicture}[scale=0.2]
	\draw[help lines] (0,0) grid (9,6);
	
	\draw [thick] (0,0) -- (1,1)--(2,2)--(3,3) --(4,4)--(5,5)--(6,6)--(7,4)--(8,2)--(9,0);
\end{tikzpicture}

	\caption{All 12 2-Dyck paths of 9 steps.}
\end{figure}

The new paper by Kariuki and Okoth \cite{Otho} offers the new idea to introduce left down-steps. If $L=(-1,-2)$, this is not very interesting, as, related to the first paper on skew Dyck paths by Deutsch et al. \cite{Emor} one should be careful that the path does not intersect (overlaps) itself, which in such a setting it wouldn t do. Thus it is the right idea to stretch down-steps
to $D=(2,-2)$ and $L=(-2,-2)$. Note that down-steps resp.\ down-left-steps still only count as one step. As we did before in our generating function based analysis of skew Dyck paths~\cite{korea}, instead of $(-2,-2)$ we draw $(2,-2)$ in red color. The syntactic conditions to be fulfilled are that $UL$ and $LU$ are forbidden. A step $L$ can only follow a $D$-step, and cannot be followed by an $U$-step; note that the very last step may be an $L$-step.

\begin{figure}[h]
	\begin{tikzpicture}[scale=0.2]
		\draw[help lines] (0,0) grid (12,5);
		
		\draw [thick] (0,0) -- (1,1)--(2,2)--(4,0) --(5,1)--(6,2)--(8,0)--(9,1)--(10,2)--(12,0);
		
	\end{tikzpicture}
	\begin{tikzpicture}[scale=0.2]
		\draw[help lines] (0,0) grid (12,5);
		
		\draw [thick] (0,0) -- (1,1)--(2,2)--(3,3) --(4,4)--(6,2)--(7,3)--(9,1)--(10,2)--(12,0);

	\end{tikzpicture}
	\begin{tikzpicture}[scale=0.2]
		\draw[help lines] (0,0) grid (12,5);
		
		\draw [thick] (0,0) -- (1,1)--(2,2)--(3,3) --(5,1)--(6,2)--(8,0)--(9,1)--(10,2)--(12,0);

	\end{tikzpicture}
	\begin{tikzpicture}[scale=0.2]
		\draw[help lines] (0,0) grid (12,5);
		
		\draw [thick] (0,0) -- (1,1)--(2,2)--(3,3) --(4,4)--(5,5)--(7,3)--(9,1)--(10,2)--(12,0);
	\end{tikzpicture}
	\begin{tikzpicture}[scale=0.2]
		\draw[help lines] (0,0) grid (12,5);
		
		\draw [thick] (0,0) -- (1,1)--(2,2)--(3,3) --(5,1)--(6,2)--(7,3)--(8,4)--(10,2)--(12,0);
	\end{tikzpicture}
	\begin{tikzpicture}[scale=0.2]
		\draw[help lines] (0,0) grid (12,5);
		
		\draw [thick] (0,0) -- (1,1)--(2,2)--(4,0) --(5,1)--(6,2)--(7,3)--(9,1)--(10,2)--(12,0);
	\end{tikzpicture}
	
	\begin{tikzpicture}[scale=0.2]
		\draw[help lines] (0,0) grid (12,6);
		
		\draw [thick] (0,0) -- (1,1)--(2,2)--(3,3) --(4,4)--(5,5)--(7,3)--(8,4)--(10,2)--(12,0);
	\end{tikzpicture}
	\begin{tikzpicture}[scale=0.2]
		\draw[help lines] (0,0) grid (12,6);
		
		\draw [thick] (0,0) -- (1,1)--(2,2)--(3,3) --(5,1)--(6,2)--(7,3)--(9,1)--(10,2)--(12,0);
	\end{tikzpicture}
	\begin{tikzpicture}[scale=0.2]
		\draw[help lines] (0,0) grid (12,6);
		
		\draw [thick] (0,0) -- (1,1)--(2,2)--(4,0) --(5,1)--(6,2)--(7,3)--(8,4)--(10,2)--(12,0);
	\end{tikzpicture}
	\begin{tikzpicture}[scale=0.2]
		\draw[help lines] (0,0) grid (12,6);
		
		\draw [thick] (0,0) -- (1,1)--(2,2)--(3,3) --(4,4)--(6,2)--(7,3)--(8,4)--(10,2)--(12,0);
	\end{tikzpicture}
	\begin{tikzpicture}[scale=0.2]
		\draw[help lines] (0,0) grid (12,6);
		
		\draw [thick] (0,0) -- (1,1)--(2,2)--(3,3) --(4,4)--(6,2)--(8,0)--(9,1)--(10,2)--(12,0);
	\end{tikzpicture}
	\begin{tikzpicture}[scale=0.2]
		\draw[help lines] (0,0) grid (12,6);
		
		\draw [thick] (0,0) -- (1,1)--(2,2)--(3,3) --(4,4)--(5,5)--(6,6)--(8,4)--(10,2)--(12,0);
	\end{tikzpicture}

	\caption{All 12 2-Dyck paths of 9 steps, down-steps stretched.}
\end{figure}
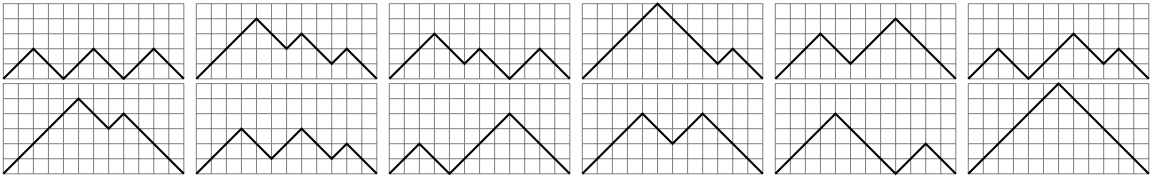

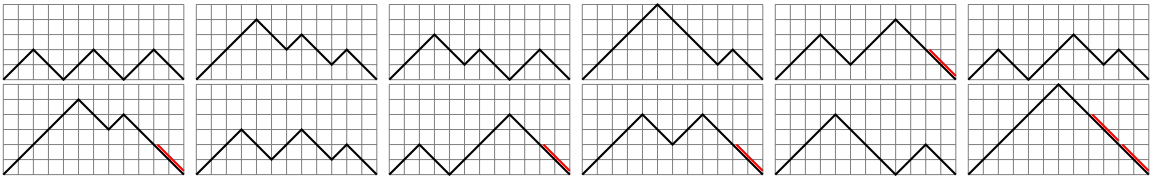
\begin{figure}[h]
	\begin{tikzpicture}[scale=0.2]
		\draw[help lines] (0,0) grid (12,5);
		
		\draw [thick] (0,0) -- (1,1)--(2,2)--(4,0) --(5,1)--(6,2)--(8,0)--(9,1)--(10,2)--(12,0);
		
	\end{tikzpicture}
	\begin{tikzpicture}[scale=0.2]
		\draw[help lines] (0,0) grid (12,5);
		
		\draw [thick] (0,0) -- (1,1)--(2,2)--(3,3) --(4,4)--(6,2)--(7,3)--(9,1)--(10,2)--(12,0);

	\end{tikzpicture}
	\begin{tikzpicture}[scale=0.2]
		\draw[help lines] (0,0) grid (12,5);
		
		\draw [thick] (0,0) -- (1,1)--(2,2)--(3,3) --(5,1)--(6,2)--(8,0)--(9,1)--(10,2)--(12,0);

	\end{tikzpicture}
	\begin{tikzpicture}[scale=0.2]
		\draw[help lines] (0,0) grid (12,5);
		
		\draw [thick] (0,0) -- (1,1)--(2,2)--(3,3) --(4,4)--(5,5)--(7,3)--(9,1)--(10,2)--(12,0);
	\end{tikzpicture}
	\begin{tikzpicture}[scale=0.2]
		\draw[help lines] (0,0) grid (12,5);
		
		\draw [thick] (0,0) -- (1,1)--(2,2)--(3,3) --(5,1)--(6,2)--(7,3)--(8,4)--(10,2)--(12,0);
					\draw [thick,red](10+0.25,2)--(12,0+0.25);
					
	\end{tikzpicture}
	\begin{tikzpicture}[scale=0.2]
		\draw[help lines] (0,0) grid (12,5);
		
		\draw [thick] (0,0) -- (1,1)--(2,2)--(4,0) --(5,1)--(6,2)--(7,3)--(9,1)--(10,2)--(12,0);
	\end{tikzpicture}
	
	\begin{tikzpicture}[scale=0.2]
		\draw[help lines] (0,0) grid (12,6);
		
		\draw [thick] (0,0) -- (1,1)--(2,2)--(3,3) --(4,4)--(5,5)--(7,3)--(8,4)--(10,2)--(12,0);
								\draw [thick,red](10+0.25,2)--(12,0+0.25);
								
	\end{tikzpicture}
	\begin{tikzpicture}[scale=0.2]
		\draw[help lines] (0,0) grid (12,6);
		
		\draw [thick] (0,0) -- (1,1)--(2,2)--(3,3) --(5,1)--(6,2)--(7,3)--(9,1)--(10,2)--(12,0);
	\end{tikzpicture}
	\begin{tikzpicture}[scale=0.2]
		\draw[help lines] (0,0) grid (12,6);
		
		\draw [thick] (0,0) -- (1,1)--(2,2)--(4,0) --(5,1)--(6,2)--(7,3)--(8,4)--(10,2)--(12,0);
				\draw [thick,red](10+0.25,2)--(12,0+0.25);
	\end{tikzpicture}
	\begin{tikzpicture}[scale=0.2]
		\draw[help lines] (0,0) grid (12,6);
		
		\draw [thick] (0,0) -- (1,1)--(2,2)--(3,3) --(4,4)--(6,2)--(7,3)--(8,4)--(10,2)--(12,0);
								\draw [thick,red](10+0.25,2)--(12,0+0.25);
								
	\end{tikzpicture}
	\begin{tikzpicture}[scale=0.2]
		\draw[help lines] (0,0) grid (12,6);
		
		\draw [thick] (0,0) -- (1,1)--(2,2)--(3,3) --(4,4)--(6,2)--(8,0)--(9,1)--(10,2)--(12,0);
	\end{tikzpicture}
	\begin{tikzpicture}[scale=0.2]
		\draw[help lines] (0,0) grid (12,6);
		
		\draw [thick] (0,0) -- (1,1)--(2,2)--(3,3) --(4,4)--(5,5)--(6,6)--(8,4)--(10,2)--(12,0);
				\draw [thick,red](10+0.25,2)--(12,0+0.25);
						\draw [thick,red](8+0.25,4)--(10,2+0.25);
	\end{tikzpicture}

	\caption{All 2-Dyck paths of 9 steps, with possible red down steps indicated. This leads to altogether 19 paths.}
	\label{bursch}
\end{figure}

The rest of this paper is devoted to an analysis (=enumeration) of skew 2-Dyck paths, using the kernel method. Some
recent applications of these techniques can be found in \cite{Prodinger-kernel}, \cite{garden}, \cite{korea}.
There are some advantages: The analysis is in a sense straight forward; one does not need to invent \emph{clever} decompositions
of the paths in question, and one gets the enumeration of partial skew 2-Dyck paths as a bonus, as we have  all three generating functions
leading to level $i$ by an up-step, down-step, red down-step in an explicit form. 

In a further section, we indicate how to extend the concept from skew 2-Dyck paths to skew $t$-Dyck paths, again handling prefixes as a bonus.

In a last section, 2-Dycks are scanned from right to left. If they return to $x$-axis, the enumeration is of course the same as before.

\section{Automaton to syntactically check 2-Dyck paths with left (red) steps }

We start by constructing an automaton that checks whether the syntactic conditions are still fulfilled, i. e., $UL$ and $LU$ are forbidden. 
According to the three types of steps, we have three layers of states (current level).

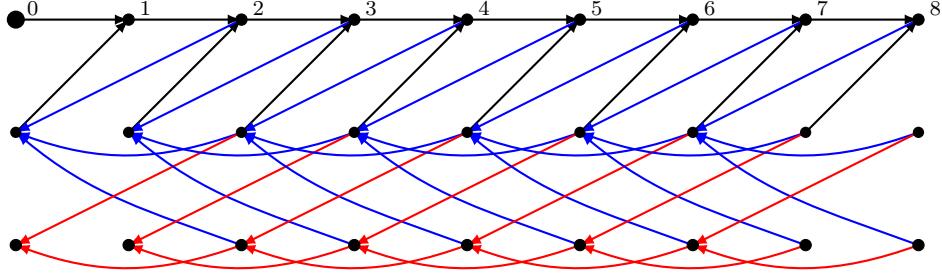
\begin{figure}[h]

	\begin{center}
		\begin{tikzpicture}[scale=1.5]

			\fill (0,0) circle (0.08cm);

			\foreach \x in {0,1,2,3,4,5,6,7}
			{
				\draw[thick, -latex] (\x,0-1) to (\x+1,0); 
				\draw[ thick, latex-] (\x+1,0) to (\x,0); 
				\node at (\x+0.15,0.1){\tiny$\x$};
			} 
			
			\foreach \x in {8}
			{
			}
			
			\foreach \x in {0,1,2,3,4,5,6}
			{
				
				\draw[ thick,blue, -latex] (\x+2,0)to (\x,-1); 
				\draw[ thick,blue, -latex] (\x+2,-1)[out=200,in=-20]to (\x,-1); 
				\draw[ thick,blue, -latex] (\x+2,-2)[out=160,in=-40]to (\x,-1); 				

				\draw[ thick,red, -latex] (\x+2,-2)[in=-20,out=200]to (\x,-2); 				
								\draw[ thick,red, -latex] (\x+2,-1)to (\x,-2); 				
			} 
			
			\node at (8+0.15,0.1){\tiny$8$};
			
			\foreach \x in {0,1,2,3,4,5,6,7,8}
			{
				\draw (\x,0) circle (0.05cm);
				\fill (\x,0) circle (0.05cm);
				\
				draw (\x,-1) circle (0.05cm);
				\fill (\x,-1) circle (0.05cm);
				\draw (\x,-2) circle (0.05cm);
				\fill (\x,-2) circle (0.05cm);
			}
		\end{tikzpicture}
	\end{center}
	\caption{Graph (automaton) to recognize skew 2-Dyck paths.}
	\label{purpelUDDU}
\end{figure}
For the enumeration, we introduce $f_i=f_i(z)$, $g_i=g_i(z)$, $h_i=h_i(z)$, the generating functions of skew 2-Dyck paths leading to state $i$ in the first (second, third layer, where the exponent $z^k$ counts the number of steps when starting at the origin (thick node on the first layer). The recursions are as follows:
\begin{align*}
f_0&=1, \quad f_{i+1}=zf_i+zg_i,\ i\ge0,\\*
g_i&=zf_{i+2}+zg_{i+2}+zh_{i+2},\ i\ge0,\\*
h_i&=zg_{i+2}+zh_{i+2},\ i\ge0.
\end{align*}
They are easy to understand by distinguishing the last possible step that leads to level $i$. 
Further, one introduces bivariate generating functions
\begin{equation*}
F(u)=F(u,z)=\sum_{i\ge0}u^if_i,\quad
G(u)=G(u,z)=\sum_{i\ge0}u^ig_i,\quad
H(u)=H(u,z)=\sum_{i\ge0}u^ih_i.
\end{equation*}
By summing the recursions, we find
\begin{align*}
F-1&=uz(F+G),\\
u^2G&=z(F-1-uf_1)+z(G-g_0-ug_1)+z(H-h_0-uh_1),\\
u^2H&=z(G-g_0-ug_1)+z(H-h_0-uh_1).
\end{align*}
The solution is (obtained with the help of a computer)
\begin{align*}
F&=\frac{z^2f_1u^3+z^2g_1u^3+z^2h_1u^3-u^3+z^2u^2+z^2g_0u^2+z^2h_0u^2-z^3f_1u+2zu-z^3}{zu^4-u^3-z^2u^2+2zu-z^3},\\
G&=\frac{z\Bigl(\genfrac{}{}{0pt}{}{h_1u^3+zf_1u^3+zg_1u^3-h_1u^2-g_1u^2+zu^2+zh_0u^2}{-f_1u^2+zg_0u^2-h_0u-g_0u-z^2f_1u+zf_1-z^2}\Bigr)}{zu^4-u^3-z^2u^2+2zu-z^3},\\
H&=\frac{z\Bigl(\genfrac{}{}{0pt}{}{ zh_1u^3+zg_1u^3-g_1u^2+zg_0u^2-h_1u^2+zh_0u^2+z^2h_1u}{+z^2g_1u+z^2f_1u-g_0u-h_0u+z^2h_0+z^2g_0-zf_1+z^2}\Bigr)}{zu^4-u^3-z^2u^2+2zu-z^3}.
\end{align*}
First, we need to understand the denominator:
\begin{equation*}
	zu^4-u^3-z^2u^2+2zu-z^3=z(u-s_1)(u-s_2)(u-s_3)(u-s_4),
\end{equation*}
with
\begin{align*}
	s_1&={\frac {1}{2}}{z}^{2}+{\frac {3}{16}}{z}^{5}+{\frac {17}{128}}{z}^{8}
	+{\frac {29}{256}}{z}^{11}+{\frac {861}{8192}}{z}^{14}+{\frac {6675}{
			65536}}{z}^{17}+{\frac {13231}{131072}}{z}^{20}+{\frac {52939}{524288}
	}{z}^{23}+\dots,\\*
s_2&=\sqrt {2}\sqrt {z}+\frac14{z}^{2}+{\frac {21}{64}}\sqrt {2}{z}^{7/2}+{
	\frac {29}{32}}{z}^{5}+{\frac {10727}{8192}}\sqrt {2}{z}^{13/2}+{
	\frac {1007}{256}}{z}^{8}+\dots,\\*
s_3&=-\sqrt {2}\sqrt {z}+\frac14{z}^{2}-{\frac {21}{64}}\sqrt {2}{z}^{7/2}+
{\frac {29}{32}}{z}^{5}-{\frac {10727}{8192}}\sqrt {2}{z}^{13/2}+{
	\frac {1007}{256}}{z}^{8}-\dots\\*
s_4&=\frac1z-{z}^{2}-2{z}^{5}-8{z}^{8}-39{z}^{11}-210{z}^{14}-
1203{z}^{17}-7192{z}^{20}-44362{z}^{23}-280250{z}^{26}-\dots
	\end{align*}
These solutions where found using \textsf{gfun} and its \textsf{algeqtoseries} construction.\footnote{Thanks to Cyril Banderier for valuable hints.}

None of the factors $u-s_1$, $u-s_2$, $u-s_3$ can appear since we expect power series solutions, so they must cancel out. Hence, in $F$, $G$, $H$ the factor
$(u-s_1)(u-s_2)(u-s_3)$ can be divided out from both, numerator and denominator,\footnote{Such a reasoning belongs to the realm of the kernel method.} with the result 
\begin{align*}
F(u)&=\frac{z^2(f_1+g_1+h_1)-1}{z(u-s_4)},\\*
G(u)&=\frac{-z^2(f_1+g_1+h_1)}{z(u-s_4)},\\*
H(u)&=\frac{-z^2(g_1+h_1)}{z(u-s_4)}.
\end{align*}
Now one can plug in $u=0$, to find
\begin{align*}
	1=f_0&=\frac{z^2(f_1+g_1+h_1)-1}{z(-s_4)},\\
	g_0&=\frac{-z^2(f_1+g_1+h_1)}{z(-s_4)},\\
	h_0&=\frac{-z^2(g_1+h_1)}{z(-s_4)}.
\end{align*}
Simplifying the system leads to 
\begin{equation*}
1+g_0= \frac1{zs_4},\quad g_0-h_0=\frac{z^2(1+g_0)}{s_4}
\end{equation*}
and further to
\begin{align*}
g_0&=\frac1{zs_4}-1={z}^{3}+3{z}^{6}+13{z}^{9}+66{z}^{12}+365{z}^{15}+2131{z}^{
	18}+12921{z}^{21}+\dots,\\
h_0&=\frac1{zs_4}-\frac{z}{s_4^2}-1={z}^{6}+6{z}^{9}+34{z}^{12}+198{z}^{15}+1191{z}^{18}+7364{z
}^{21}+\dots\;.
\end{align*}
In total
\begin{equation*}
1+g_0+h_0=1+{z}^{3}+4{z}^{6}+\mathbf{19}{z}^{9}+100{z}^{12}+563{z}^{15}+3322{z
}^{18}+20285{z}^{21}+\dots\;.
\end{equation*}
The last series is the generating function of \emph{all} skew 2-Dyck paths, and the number 19 refers to all such paths of length 9 as given in the examples. 

Note that
\begin{equation*}
g_1+h_1=\frac{h_0s_4}{z}, \quad f_1=z+zg_0,
\end{equation*}
which was used in the simplifications. Now, we can continue:
\begin{equation*}
F(u)=\frac{z^2f_1+z^2g_1+z^2h_1-1}{z(u-s_4)}=\frac{z^3(1+g_0)+zh_0s_4-1}{z(-s_4)}\frac1{1-u/s_4},
\end{equation*}
\begin{equation*}
G(u)=\frac{-z(f_1+g_1+h_1)}{u-s_4}=\frac{-z^2(1+g_0)-h_0s_4}{u-s_4}=\frac{z^2(1+g_0)+h_0s_4}{s_4}\frac1{1-u/s_4},
\end{equation*}
\begin{equation*}
	H(u)=\frac{-z(g_1+h_1)}{u-s_4}=\frac{-h_0s_4}{u-s_4}=h_0\frac1{1-u/s_4}.
\end{equation*}
Consequently
\begin{align*}
f_k&=[u^k]F(u)=\frac{-z^3(1+g_0)-zh_0s_4+1}{z}s_4^{-k-1},\\
g_k&=[u^k]G(u)=(z^2(1+g_0)+h_0s_4)s_4^{-k-1},\\
h_k&=[u^k]H(u)=h_0s_4^{-k}.
\end{align*}

These are the promised generating functions of partial skew 2-Dyck paths, leading to level $k$ in all three possible ways.

The following comments might be useful. The denominator $zu^4-u^3-z^2u^2+2zu-z^3$ means
that
\begin{equation*}
zf_{n}-f_{n+1}-z^2f_{n+2}+2zf_{n+3}-z^3f_{n+4}=0
\end{equation*}
and the other sequences $g_n$, $h_n$ satisfy the same recursion.

Surprisingly, the numbers of skew 2-Dyck paths of length $3n$ have a simple generating function (seq A007564 in \cite{OEIS}):
\begin{equation*}
R(z)=\frac{1}{6z}+\frac13-\frac{\sqrt{1-8z+4z^2}}{6z}=1+z+4z^2+19z^3+100z^4+562z^5+3304z^6+20071z^7+\dots
\end{equation*}
The extraction of coefficients is easy when using the substitution $z=\dfrac{t}{(1+t)(1+3t)}$; then $R=1+t$, and contour integration (or the Lagrange inversion formula) leads to
\begin{align*}
[z^n]R(z)&=\frac{1}{2\pi i}\oint\frac{dz}{z^{n+1}}R\\
&=\frac{1}{2\pi i}\oint dt\frac{1-3t^2}{(1+t)^2(1+3t)^2}\frac{(1+t)^{n+1}(1+3t)^{n+1}}{t^{n+1}}(1+t)\\
&=\frac{1}{2\pi i}\oint \frac{dt}{t^{n+1}}[(1+t)-t(1+3t)]{(1+t)^{n}(1+3t)^{n-1}}\\
&=\frac{1}{2\pi i}\oint \frac{dt}{t^{n+1}}(1+t){(1+t)^{n}(1+3t)^{n-1}}-\frac{1}{2\pi i}\oint \frac{dt}{t^{n+1}}t(1+3t){(1+t)^{n}(1+3t)^{n-1}}\\
&=[t^n](1+t)^{n+1}(1+3t)^{n-1}-[t^{n-1}](1+t)^{n}(1+3t)^{n}\\
&=\sum_{i=0}^n3^i\binom{n-1}{i}\binom{n+1}{n-i}-\sum_{i=0}^{n-1}3^i\binom{n}{i}\binom{n}{n-1-i}\\
&=\sum_{i=0}^n3^i\binom{n-1}{i}\binom{n+1}{i+1}-\sum_{i=0}^{n-1}3^i\binom{n}{i}\binom{n}{i+1}\\
&=\frac1n\sum_{i=0}^{n-1}3^i\binom{n}{i}\binom{n}{i+1}.
\end{align*}
The products of binomial coefficients may be written in terms of Narayana numbers $\frac1n\binom{n}{i}\binom{n}{i+1}$, compare \cite{Otho}.

The enumeration of 2-Dyck paths (without the skew-extension) leads to a third order algebraic equation, which, however, can be handled with the Lagrange inversion formula, see e. g., the thesis of Selkirk \cite{selkirk}.

\section{Generalization: skew $t$-ary Dyck-paths}

We start with skew 3-Dyck paths, from which the general case about skew
$t$-Dyck paths will become clear almost immediately. There is the usual up-step $U$, and a down-step $D$ of 3 units, and a red down-step $L$ of 3 units; the conditions are that $UL$ and $LU$ are forbidden.
\begin{figure}[h]

	\begin{center}
		\begin{tikzpicture}[scale=1.5]

			\fill (0,0) circle (0.08cm);

			\foreach \x in {0,1,2,3,4,5,6,7}
			{
				\draw[thick, -latex] (\x,0-1) to (\x+1,0); 
				\draw[ thick, latex-] (\x+1,0) to (\x,0); 
				\node at (\x+0.15,0.1){\tiny$\x$};
			}

			\foreach \x in {0,1,2,3,4,5}
			{
				
				\draw[ thick,blue, -latex] (\x+3,0)to (\x,-1); 
				\draw[ thick,blue, -latex] (\x+3,-1)[out=195,in=-20]to (\x,-1); 
				\draw[ thick,blue, -latex] (\x+3,-2)[out=170,in=-30]to (\x,-1); 				
				
				\draw[ thick,red, -latex] (\x+3,-2)[in=-20,out=200]to (\x,-2); 				
				\draw[ thick,red, -latex] (\x+3,-1)to (\x,-2); 				
			} 
			
			\node at (8+0.15,0.1){\tiny$8$};
			
			\foreach \x in {0,1,2,3,4,5,6,7,8}
			{
				\draw (\x,0) circle (0.05cm);
				\fill (\x,0) circle (0.05cm);
				\
				draw (\x,-1) circle (0.05cm);
				\fill (\x,-1) circle (0.05cm);
				\draw (\x,-2) circle (0.05cm);
				\fill (\x,-2) circle (0.05cm);
			}
		\end{tikzpicture}
	\end{center}
	\caption{Graph (automaton) to recognize skew 3-Dyck paths.}
	\label{purpel-t}
\end{figure}
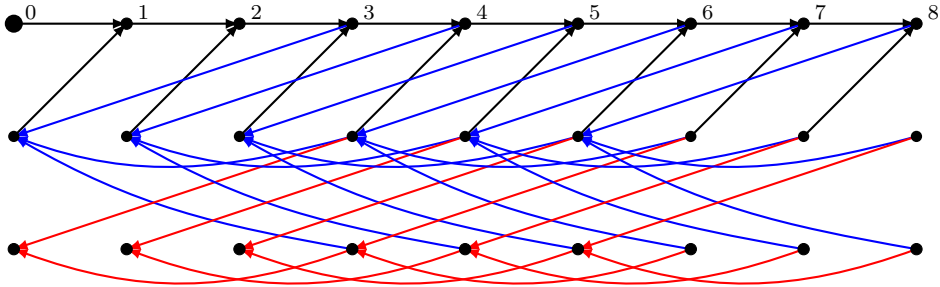

By setting up recursions and summing them, we get
\begin{align*}
	F-1&=uz(F+G),\\*
	u^3G&=z(F-1-uf_1-u^2f_2)+z(G-g_0-ug_1-u^2g_2)+z(H-h_0-uh_1-u^2h_2),\\*
	u^3H&=z(G-g_0-ug_1-u^2g_2)+z(H-h_0-uh_1-u^2h_2).
\end{align*}
The solution is
\begin{align*}
F&=\frac{u^5(z^2f_2+z^2g_2+z^2h_2-1)+\mathsf{something}}{zu^6-u^5-z^2u^3+2zu^2-z^3},\\
G&=\frac{-u^5z^2(f_2+g_2+h_2)+\mathsf{something}}{zu^6-u^5-z^2u^3+2zu^2-z^3},\\
H&=\frac{-u^5z^2(g_2+h_2)+\mathsf{something}}{zu^6-u^5-z^2u^3+2zu^2-z^3}.
\end{align*}
The denominator factors into 5 bad and one good root,
\begin{equation*}
zu^6-u^5-z^2u^3+2zu^2-z^3=z(u-s_1)(u-s_2)(u-s_3)(u-s_4)(u-s_5)(u-s_6).
\end{equation*}
The useful root is
\begin{equation*}
s_6=\frac1z-z^3-z^7-16z^{11}-104z^{15}-749z^{19}-5748z^{23}-46069z^{27}-38109z^{31}-\dots\;.
\end{equation*}
The factor $(u-s_1)(u-s_2)(u-s_3)(u-s_4)(u-s_5)$ can be divided out from numerator resp.\ denominator, with the result
\begin{align*}
F(u)&=\textsf{something}\;\frac1{1-u/s_6},\\
G(u)&=\textsf{something}\;\frac1{1-u/s_6},\\
H(u)&=\textsf{something}\;\frac1{1-u/s_6},
\end{align*}
and therefore
\begin{align*}
	f_k=[u^k]F(u)&=\textsf{something}\;s_6^{-k},\\
	g_k=[u^k]G(u)&=\textsf{something}\;s_6^{-k},\\
	h_k=[u^k]H(u)&=\textsf{something}\;s_6^{-k}.
\end{align*}

It is now clear how this translates to skew $t$-ary Dyck paths:
\begin{align*}
	f_k=[u^k]F(u)&=\textsf{something}\;s_{2t}^{-k},\\
	g_k=[u^k]G(u)&=\textsf{something}\;s_{2t}^{-k},\\
	h_k=[u^k]H(u)&=\textsf{something}\;s_{2t}^{-k};
\end{align*}
$s_{2t}$ is the only useful root of an algebraic equation of order $2t$ (in the variable $u$).

Perhaps some eager students might like to write programs so that the relevant results pop out when one enters the value of $t$.

\section{2-Dyck paths read from right to left}

First, we mirror the previous Figure \ref{bursch}. For the enumeration of 2-Dyck paths this mirrored version must lead to the same results. However, when considering partial 2-Dyck paths, results will be different.

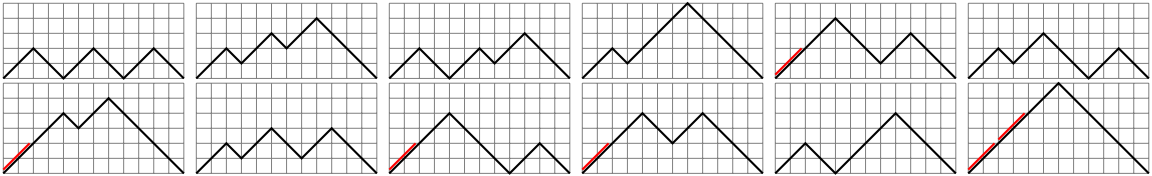
\begin{figure}[h]
	
	\begin{tikzpicture}[scale=0.2]
		\begin{scope}[xscale=-1,yscale=1]	
			\draw[help lines] (0,0) grid (12,5);
			
			\draw [thick] (0,0) -- (1,1)--(2,2)--(4,0) --(5,1)--(6,2)--(8,0)--(9,1)--(10,2)--(12,0);
		\end{scope}
	\end{tikzpicture}
	\begin{tikzpicture}[scale=0.2]
		\begin{scope}[xscale=-1,yscale=1]	
			\draw[help lines] (0,0) grid (12,5);
			
			\draw [thick] (0,0) -- (1,1)--(2,2)--(3,3) --(4,4)--(6,2)--(7,3)--(9,1)--(10,2)--(12,0);
			
		\end{scope}
	\end{tikzpicture}
	\begin{tikzpicture}[scale=0.2]
		\begin{scope}[xscale=-1,yscale=1]	
			\draw[help lines] (0,0) grid (12,5);
			
			\draw [thick] (0,0) -- (1,1)--(2,2)--(3,3) --(5,1)--(6,2)--(8,0)--(9,1)--(10,2)--(12,0);
			
		\end{scope}
	\end{tikzpicture}
	\begin{tikzpicture}[scale=0.2]
		\begin{scope}[xscale=-1,yscale=1]	
			\draw[help lines] (0,0) grid (12,5);
			
			\draw [thick] (0,0) -- (1,1)--(2,2)--(3,3) --(4,4)--(5,5)--(7,3)--(9,1)--(10,2)--(12,0);
		\end{scope}
	\end{tikzpicture}
	\begin{tikzpicture}[scale=0.2]
		\begin{scope}[xscale=-1,yscale=1]	
			\draw[help lines] (0,0) grid (12,5);
			
			\draw [thick] (0,0) -- (1,1)--(2,2)--(3,3) --(5,1)--(6,2)--(7,3)--(8,4)--(10,2)--(12,0);
			\draw [thick,red](10+0.25,2)--(12,0+0.25);
		\end{scope}		
	\end{tikzpicture}
	\begin{tikzpicture}[scale=0.2]
		\begin{scope}[xscale=-1,yscale=1]	
			\draw[help lines] (0,0) grid (12,5);
			
			\draw [thick] (0,0) -- (1,1)--(2,2)--(4,0) --(5,1)--(6,2)--(7,3)--(9,1)--(10,2)--(12,0);
	\end{scope}	\end{tikzpicture}
	
	\begin{tikzpicture}[scale=0.2]
		\begin{scope}[xscale=-1,yscale=1]	
			\draw[help lines] (0,0) grid (12,6);
			
			\draw [thick] (0,0) -- (1,1)--(2,2)--(3,3) --(4,4)--(5,5)--(7,3)--(8,4)--(10,2)--(12,0);
			\draw [thick,red](10+0.25,2)--(12,0+0.25);
		\end{scope}
	\end{tikzpicture}
	\begin{tikzpicture}[scale=0.2]
		\begin{scope}[xscale=-1,yscale=1]	
			\draw[help lines] (0,0) grid (12,6);
			
			\draw [thick] (0,0) -- (1,1)--(2,2)--(3,3) --(5,1)--(6,2)--(7,3)--(9,1)--(10,2)--(12,0);
		\end{scope}
	\end{tikzpicture}
	\begin{tikzpicture}[scale=0.2]
		\begin{scope}[xscale=-1,yscale=1]	
			\draw[help lines] (0,0) grid (12,6);
			
			\draw [thick] (0,0) -- (1,1)--(2,2)--(4,0) --(5,1)--(6,2)--(7,3)--(8,4)--(10,2)--(12,0);
			\draw [thick,red](10+0.25,2)--(12,0+0.25);
		\end{scope}
	\end{tikzpicture}
	\begin{tikzpicture}[scale=0.2]
		\begin{scope}[xscale=-1,yscale=1]	
			\draw[help lines] (0,0) grid (12,6);
			
			\draw [thick] (0,0) -- (1,1)--(2,2)--(3,3) --(4,4)--(6,2)--(7,3)--(8,4)--(10,2)--(12,0);
			\draw [thick,red](10+0.25,2)--(12,0+0.25);
		\end{scope}
	\end{tikzpicture}
	\begin{tikzpicture}[scale=0.2]
		\begin{scope}[xscale=-1,yscale=1]	
			\draw[help lines] (0,0) grid (12,6);
			
			\draw [thick] (0,0) -- (1,1)--(2,2)--(3,3) --(4,4)--(6,2)--(8,0)--(9,1)--(10,2)--(12,0);
		\end{scope}
	\end{tikzpicture}
	\begin{tikzpicture}[scale=0.2]
		\begin{scope}[xscale=-1,yscale=1]	
			\draw[help lines] (0,0) grid (12,6);
			
			\draw [thick] (0,0) -- (1,1)--(2,2)--(3,3) --(4,4)--(5,5)--(6,6)--(8,4)--(10,2)--(12,0);
			\draw [thick,red](10+0.25,2)--(12,0+0.25);
			\draw [thick,red](8+0.25,4)--(10,2+0.25);
		\end{scope}
	\end{tikzpicture}

	\caption{All 2-Dyck paths of 9 steps, with possible red down steps indicated. Right-left version.}
\end{figure}

The previous automaton can be used, but all arrows must be reversed.

\begin{figure}[h]

	\begin{center}
		\begin{tikzpicture}[scale=1.5]

			\fill (0,-1) circle (0.08cm);
				\fill (0,-2) circle (0.08cm);

			\foreach \x in {0,1,2,3,4,5,6,7}
			{
				\draw[thick, latex-] (\x,0-1) to (\x+1,0); 
				\draw[ thick, -latex] (\x+1,0) to (\x,0); 
				\node at (\x+0.15,0.1){\tiny$\x$};
			} 
			
			\foreach \x in {8}
			{
			}
			
			\foreach \x in {0,1,2,3,4,5,6}
			{
				
				\draw[ thick,blue, latex-] (\x+2,0)to (\x,-1); 
				\draw[ thick,blue, latex-] (\x+2,-1)[out=200,in=-20]to (\x,-1); 
				\draw[ thick,blue, latex-] (\x+2,-2)[out=160,in=-40]to (\x,-1); 				
				
				\draw[ thick,red, latex-] (\x+2,-2)[in=-20,out=200]to (\x,-2); 				
				\draw[ thick,red, latex-] (\x+2,-1)to (\x,-2); 				
			} 
			
			\node at (8+0.15,0.1){\tiny$8$};
			
			\foreach \x in {0,1,2,3,4,5,6,7,8}
			{
				\draw (\x,0) circle (0.05cm);
				\fill (\x,0) circle (0.05cm);
				\
				draw (\x,-1) circle (0.05cm);
				\fill (\x,-1) circle (0.05cm);
				\draw (\x,-2) circle (0.05cm);
				\fill (\x,-2) circle (0.05cm);
			}
		\end{tikzpicture}
	\end{center}
	\caption{Graph (automaton) to recognize skew 2-Dyck paths. The first step of a non-empty path must be an up-step or a red up-step. Start can be in the second or third layer; the empty path can only be counted once.}
	\label{purpelUDDU-RL}
\end{figure}
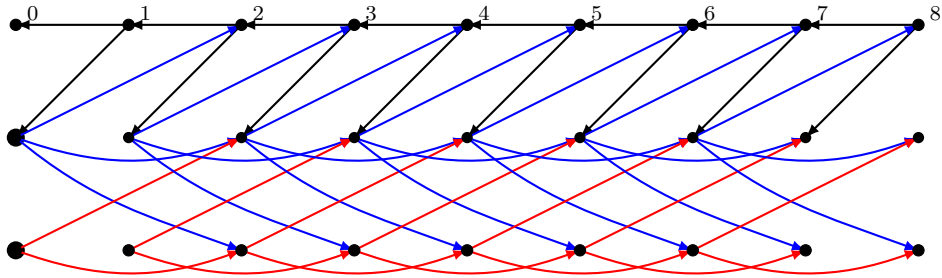
We set up the usual recursions, use the same letters, but now with a different meaning.
\begin{align*}
	f_0&=0,\quad f_1=zf_2,\quad f_{i+2}=zg_i+zf_{i+3},\ i\ge0,\\*
	g_0&=1+zf_{1},\quad g_1=zf_2,\quad g_{i+2}=zf_{i+3}+zg_i+zh_i,\ i\ge0,\\*
	h_0&=1,\quad h_1=0,\quad h_{i+2}=zg_i+zh_i,\ i\ge0.
\end{align*}
Summing,
\begin{align*}
\sum_{i\ge0}u^{i+3}f_{i+2}&=\sum_{i\ge0}u^{i+3}zg_i+\sum_{i\ge0}u^{i+3}zf_{i+3},\\*
\sum_{i\ge0}u^{i+3}g_{i+2}&=\sum_{i\ge0}u^{i+3}zf_{i+3}+\sum_{i\ge0}u^{i+3}zg_i+\sum_{i\ge0}u^{i+3}zh_i,\\*
\sum_{i\ge0}u^{i+2}h_{i+2}&=\sum_{i\ge0}u^{i+2}zg_i+\sum_{i\ge0}u^{i+2}zh_i,
\end{align*}
or
\begin{align*}
	u(F-uf_1)&=u^3zG+z(F-uf_1-u^2f_2),\\
	u(G-g_0-ug_1)&=z(F-uf_1-u^2f_2)+u^3zG+u^3zH,\\
	H-1&=u^2zG+u^2zH.
\end{align*}
The solution is
\begin{align*}
F&=\frac {u (2zu^2g_0-zu^2-g_0+1) }{-2{u}^{3}z+{z}^{2}{u}^{2}+{z}^{3}{u}^{4}+u-z},\\
	G&=\frac {-z^2u^2+u-zg_0+z^2u^2g_0}{-2{u}^{3}z+{z}^{2}{u}^
		{2}+{z}^{3}{u}^{4}+u-z},\\
	H&=\frac{u-z-z^2u^2g_0}{-2{u}^{3}z+{z}^			{2}{u}^{2}+{z}^{3}{u}^{4}+u-z}.
	\end{align*}

The factorization of the denominator  $-2u^3z+z^2u^2+z^3u^4+u-z$ is easy since it is the reciprocal denominator from before, where $u\leftrightarrow u^{-1}$. Thus we have
\begin{equation*}
	-2u^3z+z^2u^2+z^3u^4+u-z=z^3(u-s_1^{-1})(u-s_2^{-1})(u-s_3^{-1})(u-s_4^{-1}).
\end{equation*}
The only `bad' factor is $(u-t_1)$ with $t_1=s_1^{-1}$. Dividing it out from the numerator of $F$, we get $-g_0+1+2t_1^2zg_0-t_1^2z$, but since $f_0=0$, we find finally
\begin{equation*}
g_0=\frac{1-zt_1^2}{1-2zt_1^2}=-1-zt_1^2+\frac{2t_1}{z}=1+z^3+4 z^6+19 z^9+100 z^{12}+563 z^{15}+3322 z^{18}+20285 z^{21}+\dots\!.
\end{equation*}
The expressions $[u^k]G(u)$ can be computed as well, but are not as pleasant as before since we have to deal with a factor $\dfrac1{(u-s_2^{-1})(u-s_3^{-1})(u-s_4^{-1})}$,
and extraction of $u^k$ is not nice, although doable when applying partial fraction decomposition first. We leave details to interested readers.


\clearpage

\begin{thebibliography}{1}
	
	
	
	
	
	
	
	\bibitem{BF}
	C. Banderier and P. Flajolet.
	\newblock {Basic analytic combinatorics of directed lattice paths.}
	\newblock {\em Theoret. Comput. Sci.} 281 (2002), 37--80.
	
	
	\bibitem{Emor}
E. Deutsch and E. Munarini and S. Rinaldi.
	\newblock {Skew Dyck paths.}
	\newblock {\em Journal of Statistical Planning and Inference} 140 (2010), 2191--2203.
 	
	
	
 
	
	\bibitem{korea}
	H. Prodinger.
	\newblock 	Partial Skew Dyck Paths: A kernel method approch.
	\newblock {\em Graphs and Combinatorics} 38 (2022): 135 (11 pages).
	
	\bibitem{Prodinger-kernel}
	H. Prodinger.
	\newblock The kernel method: a collection of examples.
	\newblock {\em S\'{e}m. Lothar. Combin.} (2004), Art. B50f.
	
	\bibitem{garden}
	H. Prodinger.
	\newblock A walk in my lattice path garden.
	\newblock {\em S\'{e}m. Lothar. Combin.}, 87B (2023), Art. 1.
	
	\bibitem{OEIS} 
	The OEIS Foundation Inc.,
	\emph{The On-Line Encyclopedia of Integer Sequences},
	\tt{https://oeis.org} \rm (2021).
	
	
	
	\bibitem{Otho}
	Y. W. Kariuki and I. O. Okoth.
	\newblock On a class of skew {D}yck-paths.
	\newblock {\em Journal of Discrete Mathematics and Its Applications } 10 (2025), 305--319.
	
	\bibitem{selkirk}
S. J. Selkirk.
	\newblock On a generalisation of $k$-{D}yck-paths.
	\newblock {\em M. Sc.  thesis Stellenbosch} (2019).
	
	
\end{thebibliography}
\end{document}